\documentclass[smallextended]{svjour3}
\usepackage{graphicx} 
\usepackage{multirow}
\usepackage[T1]{fontenc}
\usepackage{siunitx}
\usepackage{amsmath,amsfonts,mathtools}
\usepackage{amssymb}
\usepackage{multirow}
\usepackage{hhline}
\usepackage{colortbl}
\usepackage{diagbox} 
\usepackage{makecell} 
\usepackage{slashbox}
\usepackage{subcaption}
\usepackage{derivative}
\usepackage{algorithm}
\usepackage{algpseudocode}
\usepackage{csvsimple}
\usepackage{algcompatible}

\newcommand{\mc}{\mathrm{i}}
\DeclareMathOperator{\diag}{diag}

\DeclareMathOperator{\giv}{giv}
\DeclareMathOperator{\hypf}{hypf}

\usepackage{listings}
\usepackage{subdepth} 
\newcommand{\B}[1]{\mbox{\boldmath $#1$}}
\title{Solving Shifted Systems for  Quasiseparable Matrices}
\author{Luca Gemignani}
\institute{L. Gemignani \at
              Dipartimento di Informatica, Universit\`a di Pisa \\
              Largo B. Pontecorvo 3, 56127, Pisa, Italy
              \email{luca.gemignani@unipi.it}           
}

\date{}

\begin{document}

\maketitle
\begin{abstract}
    In this paper we develop fast numerical algorithms for solving shifted linear systems  with semidefinite quasiseparable matrices.
    A combination of Givens and hyperbolic plane rotations is used to update the  Cholesky-type factorization  of the input  quasiseparable matrix by 
    determining  a factorization of its  shifted version  of the   form $LDL^T$, where $L$ is lower triangular and $D$ is a signature matrix.  
    If the shifted matrix is also definite then the  Cholesky factorization of the shifted matrix is computed in a stable way by using orthogonal transformations.  Since  quasiseparability is maintained under diagonal shifting,  a fast variant of the updating procedure using 
    computations with  generators   is  also devised.  Numerical experiments show the effectiveness and  robustness of the proposed algorithm. 
\end{abstract}

\section{Introduction}
The  solution of shifted linear systems of the form $(T^2 + \alpha I_n )\B x=\B b $, where $T\in \mathbb R^{n\times n}$ is a symmetric quasiseparable matrix,  is  at the core of the methods  devised in \cite{Aceto2025ER,Boito,Boito1,Boito2025} for solving certain  nonlocal boundary value problems. More generally, sequences of   shifted systems with   symmetric banded/quasiseparable  matrices  occur in the numerical solution of discretized nonlinear systems of ordinary and partial differential
equations by means of  implicit methods (see, e.g.,  \cite{alma990002123610203299} and \cite{Benzi2003,Meurant}  and the references therein).   Shifted linear systems are also encountered  in the application of the bisection  method for eigenvalue computation of Hermitian quasiseparable matrices \cite{etna_vol59_pp60-88}.

Several different authors have devised fast update and downdating algorithms to modify the triangular factorization of 
 displacement structured matrices by using  generalized unitary and hyperbolic transformations \cite{ALEXANDER19883,Chun,Cybenko,BOJANCZYK2000183}. In this paper we consider  extensions of these algorithms for the efficient construction of a Cholesky-like  factorization  for
shifted matrices of the form $A_\alpha=A + \alpha I_n=L_\alpha D_\alpha L_\alpha^T$
where A is a symmetric quasiseparable matrix   of order $n$, $L_\alpha$ is lower triangular, $D_\alpha$ is a signature matrix  and $\alpha$  is a
real shift.  We assume that a  Cholesky-like  factorization is initially computed for the
matrix $A$, or possibly for  $A_\alpha$ for some initial value of $\alpha$.  Observe that for $A=T^2$,  with  $T$ symmetric,  a Cholesky decomposition of $A$  can be derived from the QR factorization of $T$ without explicitly computing the squared matrix.  The question is
then how to efficiently update this factorization for subsequent shift values by exploiting the 
quasiseparable structure of both $A$ and $A_\alpha$. Moreover, it should be easier to  modify an existing  Cholesky-like factorization  in an accurate and fast way than to recompute a decomposition from scratch.

Our contribution is the design of an  updating  scheme based on using Givens  and/or hyperbolic rotations to solve the linear system $A_\alpha \B x=\B b$, without forming $A_\alpha$ explicitly.  Specifically, based on \cite{BOJANCZYK2000183}
we  begin by elaborating on the factorization 
\[
A_\alpha=\left[\begin{array}{c|c}
L_0 & \sqrt{|\alpha|} \end{array}
\right] \left[\begin{array}{cc}
D_0 \\ & \pm I_n\end{array}
\right] \left[\begin{array}{cc}
L_0^T\\ \sqrt{|\alpha|} \end{array}
\right], 
\]
where the middle factor is still a signature matrix.  Then, we  find a matrix $H\in \mathbb R^{2n \times 2n}$ such that 
\[
H \left[\begin{array}{cc}
L_0^T\\ \sqrt{|\alpha|} \end{array}
\right]=\left[\begin{array}{cc}
0\\ L_\alpha^T \end{array}
\right]
\]
with 
\[
H \left[\begin{array}{cc}
D_0 \\ & \pm I_n\end{array}
\right] H^T=\widehat D_\alpha, \quad D_\alpha=\widehat D_\alpha(n+1\colon 2n, n+1\colon 2n).
\]
In the generic case,  the matrix $H$ can   be  determined as a product of plane rotations applied to progressively annihilate the entries of $L_0^T$ according to the signature of the plane coordinates.  

Cholesky-like factorizations of quasiseparable matrices inherit the rank structure of the input matrix \cite{Yuli_book}.  If $A$ is banded  with 
bandwidth $\kappa$  then the complexity of the  updating algorithm is $O(n \kappa^2)$ flops which is  optimal with respect to both $\kappa$ and $n$.  If $A$ is quasiseparable with the quasiseparability order $\kappa$ then $L_0$ shares the same property.  Given the representation of $L_0$ via generators,  our updating procedure can be adjusted to compute  in linear time a representation of $L_\alpha$ in out-of-band quasiseparable form \cite{EIDELMAN2008266}. For a semiseparable  matrix of size $n$ and order  
$\kappa$  the updating procedure can  also be carried out at the cost of $O(n \kappa^2)$ flops. For a general quasiseparable matrix the computation of the modified  factorization requires $O(n \kappa^3)$ flops but, since the generators of the the out-of-band part  of  $L_\alpha$ basically coincide with those ones for $L_0$,    the approach yields significant computational savings in terms of leading constants.  
 If   $A$ is positive semidefinite and $\alpha>0$ then the  computation of the updated factorization can be carried out in a stable way by using orthogonal transformations only. In addition, the proposed updating procedure can also be exploited to  cheaply compute approximate incomplete factorizations  to be used as preconditioners for iterative linear solvers like PCG.

The paper is organized as follows. In Section \ref{two} we recall some basic results and definitions.  Section 
\ref{three}  develops our update algorithms for both general and structured cases.  Section \ref{four} is devoted to the results of numerical experiments that confirm the efficiency and robustness of the proposed schemes. Finally, conclusions and future work are drawn in Section \ref{five}. 

\section{Preliminaries}\label{two}

Customary  numerical methods for dense matrices  make use of  unitary transformations to reduce the input matrix in some convenient form. A generalization of  the concept of unitary  matrix 
  in vector spaces  equipped with an indefinite scalar product is provided by the following \cite{GLR}.
\begin{definition}
For a   given diagonal {\em signature } matrix   $D=\diag[\pm 1]\in \mathbb R^{n \times n}$ a matrix $Q\in \mathbb C^{n\times n}$ is called D-orthogonal if $Q^T D Q= D$. 
\end{definition} 
The reduction of a vector $\B x\in \mathbb R^n$ to a scalar multiple of the first  column of $I_n$  by using a  D-orthogonal matrix $Q$ can be carried out by annihilating one entry of $\B x$ at a time in a process akin to the use of orthogonal Givens rotations.  Thus,  we can restrict ourselves to  the case $n=2$.  
 It is also important to note
that it may be necessary to rescale the problem and then scale back the answer to
avoid intermediate overflow/underflow issues.

If $D=\pm I_2$  and $\B x^T \B x\neq 0$, $\B x=[x_1,x_2]^T$, with $|x_1|\geq |x_2|$ then the  Givens plane  rotation   matrix 
\begin{equation}\label{orthG}
 \mathcal G=\giv(x_1, x_2)\colon =\displaystyle\frac{1}{\sqrt{1+\rho^2}}\left[\begin{array}{cc}
1& \rho
\\
-\rho  & 1\end{array}
\right], \quad \displaystyle\rho=\frac{x_2}{x_1}
\end{equation}
is such that 
\[
\mathcal G \left[\begin{array}{c}
x_1\\x_2\end{array}\right]=\left[\begin{array}{c}
y_1\\0\end{array}\right],  \quad y_1=\frac{|x_1|}{x_1} \parallel \B x\parallel_2,  \quad \mathcal G^T  I_2  \  \mathcal G= I_2. 
\]

When $D=\pm \diag[1, -1]$ and $\B x^T \cdot D\cdot \B x\neq 0$ with $|x_1|>|x_2|$,  then   the  hyperbolic Givens matrix $\mathcal H$ 
\begin{equation}\label{hypG}
 \mathcal H=\hypf(x_1, x_2)\colon =\displaystyle\frac{1}{\sqrt{1-\rho^2}}\left[\begin{array}{cc}
1 & -\rho
\\
-\rho & 1\end{array}
\right], \quad \displaystyle\rho=\frac{x_2}{x_1}
\end{equation}
is D-orthogonal and maps $\B x $ as desired
\[
\mathcal H  \left[\begin{array}{c}
x_1\\x_2\end{array}\right]=\left[\begin{array}{c}
y_1\\0\end{array}\right],     
\quad y_1=\frac{|x_1|}{x_1} \sqrt{|\B x^T D \B x|}, 
\quad \mathcal H^T  D  \mathcal H= D. 
\]
It is worth noting that the matrix $\mathcal H$ in \eqref{hypG} can be decomposed as follows: 
\begin{equation}\label{eigH}
\mathcal H=\frac{1}{\sqrt{2}}\left[\begin{array}{cc}
1 & 1
\\
-1 & 1\end{array}
\right]  \diag[\sqrt{\frac{1+\rho}{1-\rho}}, \sqrt{\frac{1-\rho}{1+\rho}}]\left[\begin{array}{cc}
1 & -1
\\
1& 1\end{array}
\right]\frac{1}{\sqrt{2}}
\end{equation}
where  $Q=\giv(1,1)=\displaystyle\frac{1}{\sqrt{2}}\left[\begin{array}{cc}
1 & 1
\\
-1 & 1\end{array}
\right]$ is orthogonal. The representation \eqref{eigH}  is  exploited in \cite{CHSY} to  design a stable way to implement hyperbolic transformations.  Givens and hyperbolic plane rotation matrices are extended to the complementary case $|x_1|<|x_2|$ by setting
\[
\giv(x_1,x_2)\leftarrow \giv(x_2,x_1)\mathcal J_2 \quad \hypf(x_1,x_2)\leftarrow \hypf(x_2,x_1)\mathcal J_2, 
\]
where $\mathcal J_2=\allowdisplaybreaks\left[\begin{array}{cc}0& 1
\\
1 & 0\end{array}\right]$ is the unit antidiagonal matrix of size $2$.
If $|x_1|=|x_2|\neq 0$ then the  hyperbolic transformation  $\hypf(x_1,x_2)$ is not defined. Similar to hyperbolic rotations, hyperbolic  and unified Householder reflections \cite{ALEXANDER19883,BOJANCZYK2000183} can be introduced as generalizations of classical Householder orthogonal elementary matrices.   

It is well known that orthogonalization methods  with respect to a bilinear form are related with Cholesky-like triangular factorizations of symmetric possibly indefinite matrices (compare with \cite{ROS} and the references given therein). 

\begin{definition}
Let $A\in \mathbb R^{n\times n}$ be a symmetric possibly indefinite matrix. A Cholesky-like factorization of $A$ has the form $A=L \Omega L^T$     where $L \in \mathbb R^{n\times n}$ is lower triangular and $\Omega $  is some diagonal signature matrix.
\end{definition} 

A strongly nonsingular symmetric matrix admits an essentially unique Cholesky-like factorization.  Cholesky-like factorizations of quasiseparable matrices inherit the rank structure of the input matrix.  Let us recall the following definition \cite{Yuli_book}.

\begin{definition}\label{quasi}
A matrix $A=(a_{i,j})\mathbb R^{n\times n}$ is {\em quasiseparable} of order $\kappa$ if  there exist lower generators 
$\B p_i\in \mathbb R^{\kappa}$, $2\leq i\leq n$, $\B q_j\in \mathbb R^{\kappa}$, $1\leq j\leq n-1$,
$A_\ell\in \mathbb R^{\kappa \times \kappa}$, $2\leq \ell\leq n-1$ and upper generators 
$\B g_i\in \mathbb R^{\kappa}$, $1\leq i\leq n-1$, $\B h_j\in \mathbb R^{\kappa}$, $2\leq j\leq n$,
$B_\ell\in \mathbb R^{\kappa \times \kappa}$, $2\leq \ell\leq n-1$, such that 
\[
a_{i,j}=\left\{\begin{array}{ll}
\B p_i^T \overleftarrow{A_{i,j}} \B q_j \ {\rm if} \ i>j; \\
\B g_i^T \overrightarrow{B_{i,j}} \B h_j \ {\rm if} \ i<j;
\end{array}
\right.
\]
where we denote
\[
\overleftarrow{A_{i,j}}= \left\{\begin{array}{ll}
A_{i-1}\cdots A_{j+1} \ {\rm if} \ i>j+1; \\
I_\kappa \ {\rm if} \ i=j+1;
\end{array}
\right.   \quad \overrightarrow{B_{i,j}}= \left\{\begin{array}{ll}
B_{i+1}\cdots B_{j-1} \ {\rm if} \ j>i+1; \\
I_\kappa \ {\rm if} \ j=i+1;
\end{array}
\right.  
\]
When $A_i=B_i=I_\kappa$, $2\leq i\leq n$, then $A$ is  referred to as a {\em semiseparable} matrix of order $\kappa$.
\end{definition}
If $A$ is quasiseparable of order  $\kappa$  and  there exists a Cholesky-like factorization of $A=L \Omega L^T$ with $L$ invertible, then $L$ is also quasiseparable of order less than or equal to $\kappa$.  In the next section we investigate the problem of computing a Cholesky-like triangular factorization  of $A+\alpha I_n$ given the initial decomposition of the quasiseparable matrix $A$. 

\section{A Fast Method to Update the Triangular Factorization of a Shifted Quasiseparable Matrix}\label{three}
Given  a symmetric possibly indefinite matrix  $A\in \mathbb R^{n\times n}$  let $A=A_0=L_0\Omega_0 L_0^T$ be  its  Cholesky-like triangular factorization, where $L_0$ is lower triangular and $\Omega_0$ is some signature matrix.  Here we are interested in computing a 
Cholesky-like triangular factorization of the shifted matrix $A_\alpha=A+\alpha I_n=L_\alpha \Omega_\alpha L_\alpha^T$ with $\alpha \in \mathbb R$.  Based on \cite{BOJANCZYK2000183}, we first describe the general strategy  for the modified Cholesky problem and then show how to adjust this strategy for the case of quasiseparable matrices. 

Let us start by observing that 
\[
A_\alpha=\left[\begin{array}{c|c}
L_0 & \sqrt{|\alpha|} \end{array}
\right] \left[\begin{array}{cc}
\Omega_0 \\ & \pm I\end{array}
\right] \left[\begin{array}{cc}
L_0^T\\ \sqrt{|\alpha|} \end{array}
\right]=\left[\begin{array}{c|c}
L_0 & \sqrt{|\alpha|} \end{array}
\right] \widehat \Omega_0\left[\begin{array}{cc}
L_0^T\\ \sqrt{|\alpha|} \end{array}
\right].
\]
The matrix $\widehat \Omega_0$ is still a signature matrix which makes possible to devise an orthogonalization-like process to incrementally annihilate the entries of  $L_0$ and $L_0^T$.  Specifically, let  us define $R_0=\left[\begin{array}{c|c}
L_0 & \sqrt{|\alpha|} \end{array} 
\right]^T$ such that 
\begin{equation}\label{startf}
A_\alpha=A+\alpha I_n=R_0^T \widehat \Omega_0 R_0.
\end{equation}
At the first step we annihilate the entry  $R_0(1,1)$ by means a rotation in the plane of coordinates 1 and $n+1$.   The transformation is  carried out  by a Givens matrix if $\widehat \Omega_0(1,1)$ and $ \widehat \Omega_0(n+1,n+1)$  have the same sign  or  by an hyperbolic  rotation if, otherwise, $\widehat \Omega(1,1)$ and $ \widehat \Omega(n+1,n+1)$  have different signs. By setting 
\[
Q_1=I_{2n} + \left[ \B e_1, \B e_{n+1} \right ]\mathcal U_1 \left[ \B e_1, \B e_{n+1} \right ]^T
\]
where $\mathcal U_1\in \{\mathcal G,\mathcal H,   \mathcal G\mathcal J_2, \mathcal H\mathcal J_2\}$, we find that 
\[
A_\alpha=R_0^T \widehat \Omega_0 R_0=R_0^T Q_1^T (Q_1^{-T}\widehat \Omega_0 Q_1^{-1}) Q_1 R_0= R_1^T \widehat \Omega_1 R_1. 
\]
It can be easily seen that  $\widehat \Omega_1$ is a signature matrix. Specifically, 
$\widehat \Omega_1=\widehat \Omega_0 \iff \mathcal U_1\in \{\mathcal G,\mathcal H \mathcal J_2,  \mathcal G\mathcal J_2\}$,  whereas 
$\widehat \Omega_1(1,1)= \widehat \Omega_0(n+1,n+1)$ and $ \widehat \Omega_1(n+1,n+1)=\widehat \Omega_0(1,1)$  for $\mathcal U_1=\mathcal H$. Pictorially,  for $n=4$ the matrix $R_1$ reads as follows: 
 \[
\left[\begin{array}{cccc}
 & x & x & x \\ & x& x& x \\ & & x& x\\& & & x  \\ \hline x& x& x& x \\& x \\ & & x\\ & & & x      \end{array}     
\right].
 \]
The elimination  scheme can thus be continued by performing two rotation in the plane of coordinates $(1,n+2)$ and $(2, n+2)$ to zero the entries $R_1(1,2)$ and $R_1(2,2)$, respectively. 
Proceeding in this way,  whenever the process  goes to completion, at the very end  we would find  a matrix $R_n$
with 
\[
R_n=\left[\begin{array}{c|c}
0 & L_\alpha \end{array} 
\right]^T
\]
and a  diagonal signature matrix $\widehat \Omega_n$  such that 
\[
A_\alpha=R_n^T\widehat \Omega_n 
R_n=L_\alpha\widehat \Omega_n(n+1\colon 2n, n+1\colon 2n) L_\alpha^T =L_\alpha \Omega_\alpha L_\alpha^T.
\]
The following  two pseudocodes describe the resulting algorithms. 

\begin{algorithm}[H]
 \scriptsize{
\caption{: Procedure {\textsc{ZEROING}}. Given $R\in \mathbb R^{2n\times n}$, $\Omega \in \mathbb R^{2n\times 2n}$ and two indices $1\leq k\leq j\leq n$, this  procedure annihilate the entry of $R$ in position $(k,j)$ by returning  the modified $R$ and $\Omega$.}
    \label{algorithm3}
    \begin{algorithmic}[1] 
    \STATE $\B x=[R(k,j); R(n+j,j)]$; 
       \IF {$\Omega([k,n+j], [k,n+j])=\pm I_2$}
         \STATE  Determine $\mathcal U\in \{\mathcal G, \mathcal G \mathcal J_2\}$ such that  $\mathcal U \B x=[y; 0]$; 
         \STATE $R([k,n+j], \colon)=\mathcal J_2 \mathcal UR([k,n+j], \colon)$; 
         \ELSE
         \IF {$|x_1|>|x_2|$}
          \STATE  Determine $\mathcal U=\mathcal H$ such that  $\mathcal U \B x=[y; 0]$;
          \STATE $R([k,n+j], \colon)=\mathcal J_2 \mathcal UR([k,n+j], \colon)$; 
          \STATE $\Omega([k,n+j], [k,n+j])=\mathcal J_2 \Omega([k,n+j], [k,n+j]) \mathcal J_2$; 
          \ELSE
          \STATE  Determine $\mathcal U=\mathcal H \mathcal J_2$ such that  $\mathcal U \B x=[y; 0]$;
          \STATE $R([k,n+j], \colon)=\mathcal J_2 \mathcal UR([k,n+j], \colon)$;
\ENDIF
         \ENDIF
    \end{algorithmic}}
\end{algorithm}
\vspace{-1cm}
\begin{algorithm}[H]
 \scriptsize{
\caption{: Procedure {\textsc{UPDATE}}.  Given $R=R_0\in \mathbb R^{2n\times n}$ and $\Omega =\widehat \Omega_0\in \mathbb R^{2n\times 2n}$ as defined in \eqref{startf} this  algorithm  computes a Cholesky-like factorization of the shifted matrix $A_\alpha=R^T \Omega R=L_\alpha \Omega_\alpha L_\alpha^T$.}
    \label{algorithm4}
    \begin{algorithmic}[1] 
     \FOR  {$j=1\colon n$};
         \FOR  {$k=1\colon j$};
          \STATE  $[R, \Omega]=ZEROING(R, \Omega, k, j)$; 
          \ENDFOR
          \ENDFOR
          \STATE $L_\alpha=(R(n+1\colon 2n, :))^T$, $\Omega_\alpha=\Omega(n+1\colon n, n+1\colon n)$; 
    \end{algorithmic}}
\end{algorithm}

\begin{remark}
 Concerning the applicability of  Procedure {\textsc{UPDATE}} some clarifications are in order.   The annihilation scheme based on generalized plane rotations  applied in a specified order can break down prematurely. However, by a continuity argument  it can be proved that the annihilation scheme based on generalized Householder transformations  works for any strongly nonsingular matrix $A_\alpha$.  The generalized Householder matrix can  thus be decomposed as product of plane rotations by revealing the existence of an annihilation scheme  that also uses  generalized plane rotations. In our approach we prefer to deal with  schemes based on plane rotations applied in a given specified order since it simplifies the description of the process for  banded and more generally quasiseparable matrices. Generalizations of these schemes using hyperbolic Householder matrices will be treated elsewhere.
\end{remark}

Procedure {\textsc{UPDATE}} can also be extended to compute an incomplete Cholesky-like factorization according to the  sparsity pattern of $L_0^T$. Let us define 
the following sets which specify the sparsity pattern of the matrix $L_0$:
$\Delta=\{(k,j)\colon 1\leq k\leq j\leq n\}$, $\mathcal S=\{(k,j)\in \Delta \colon (L_0^T)_{k,j}=0\}$, $\mathcal S^\complement=\{(k,j)\in \Delta \colon (L_0^T)_{k,j}\neq 0\}$. Moreover, set $\mathcal R_k=\{j\colon (k,j)\in \mathcal S\}$, 
$\mathcal C_j=\{k\colon (k,j)\in \mathcal S\}$,  $\mathcal R_k^\complement=\{j\colon (k,j)\in \mathcal S^\complement\}$ and 
$\mathcal C_j^\complement=\{k\colon (k,j)\in \mathcal S^\complement\}$.  
For a given vector $\B v=\left[0, \ldots, 0, v_k, \ldots, v_n\right]$  the row-projection operator $\mathcal P_k$ is such that 
$\mathcal P_k(\B v)=\left[0, \ldots, 0, w_k, \ldots, w_n\right]$ where 
$w_s=0$ if $(k, s)\in \mathcal R_k$ and $w_s=v_s$ if, otherwise,  $(k, s)\in \mathcal R_k^\complement$. 

The following modification of Procedure {\textsc{ZEROING}} can be used to mantain the sparsity pattern of $L_0^T$. This modified version can be incorporated in  Procedure {\textsc{UPDATE}} to obtain an incomplete updating scheme. 

\begin{algorithm}[H]
 \scriptsize{
\caption{: Procedure {\textsc{INZEROING}}. Given $R\in \mathbb R^{2n\times n}$, $\Omega \in \mathbb R^{2n\times 2n}$ and two indices $1\leq k\leq j\leq n$, $(k,j)\in \mathcal S^\complement $ this  procedure annihilate the entry of $R$ in position $(k,j)$ by returning  $\Omega$ and   the modified $R$ according to the sparsity pattern  specified by $\mathcal S$.}
    \label{algorithm5}
    \begin{algorithmic}[1] 
    \STATE  $[R, \Omega]=ZEROING(R, \Omega, k, j)$; 
          \STATE $R(k, \colon)=\mathcal P_k(R(k, \colon))$; $R(n+j, \colon)=\mathcal P_j(R(n+j, \colon))$; 
    \end{algorithmic}}
\end{algorithm}
\vspace{-1cm}
\begin{algorithm}[H]
 \scriptsize{
\caption{: Procedure {\textsc{INUPDATE}}.  Given $R=R_0\in \mathbb R^{2n\times n}$ and $\Omega =\widehat \Omega_0\in \mathbb R^{2n\times 2n}$ as defined in \eqref{startf} this  algorithm  computes an incomplete  Cholesky-like factorization of the shifted matrix $A_\alpha=R^T \Omega R=L_\alpha \Omega_\alpha L_\alpha^T$ according to the sparsity pattern $\mathcal S$ of $L_0^T$.}
    \label{algorithm6}
    \begin{algorithmic}[1] 
     \FOR  {$j=1\colon n$};
         \FOR  {$k\in \mathcal C_j^\complement$};
          \STATE  $[R, \Omega]=INZEROING(R, \Omega, k, j)$; 
          \ENDFOR
          \ENDFOR
          \STATE $L_\alpha=(R(n+1\colon 2n, :))^T$, $\Omega_\alpha=\Omega(n+1\colon n, n+1\colon n)$; 
    \end{algorithmic}}
\end{algorithm}

\vspace{-0.5cm}

The  worst-case  complexity of {\textsc{UPDATE}} is $O(n^3)$  flops  but it  can decrease significantly  under  suitable assumptions. In particular, if $R_0$ is banded  with bandwidth $\kappa$ then the complexity becomes $O(n \kappa^2)$.  Banded triangular factors   are calculated while factoring  banded matrices.  The class of quasiseparable matrices includes banded matrices as a subclass. This  fact supports the   adaptation of {\textsc{UPDATE}}   to the  more general class of quasiseparable matrices that will be discussed in the next subsection. The complexity of {\textsc{INUPDATE}} is proportional to the number of non-zeros in the matrix $L_0^T$  by approaching $O(n)$  for  strict banded  triangular factors.

\subsection{ A Fast Generator-based Reduction Scheme}
A fast algorithm for computing the updated triangular factorization of $A_\alpha=A+\alpha I$ where $A$ is quasiseparable  can be  obtained by interlacing  band reduction stages   with  the action of a bulge chasing scheme.  The reduction in banded form is accomplished by a process like  the modification of the Dewilde-van der Veen method described in \cite{EIDELMAN2002419}. In the band reduction stages  some  computations can be performed  offline and stored thus allowing for computational savings in the updating procedure.

To illustrate the  intuitive idea and for the sake of simplicity 
let us begin by  first considering  quasiseparable matrices of 
quasiseparability rank $\kappa=1$. Let $R=R_0$ be specified as in \eqref{startf} where $L_0^T$ is quasiseparable with upper quasiseparability rank $\kappa=1$.  Firstly,   we apply a plane rotation  in the plane of coordinates $(1,2)$ to zeroing  some entries in the first row of $L_0^T$. Specifically, due to the rank structure we  may reduce $R$ in the form 
\vspace{-0.1cm}
 \[
R=\left[\begin{array}{cccc}
 x& x & 0 & 0 \\ x& x& x& x \\ & & x& x\\& & & x  \\ \hline  x\\& x \\ & & x\\ & & & x      \end{array}     
\right] .
 \]
 The annihilation in the first row generates a bulge in the top-left corner of the matrix $R$.
 Then,  we   nullify the first column of this  bulge by  appling  two plane rotations in the plane $(1, n+1)$ and $(2,n+1)$ to annihilate $R(1,1)$ and $R(2,1)$, respectively.  
 It should be noticed that  the  newly generated row in the botton part of $R$  is a linear combination of the first two rows and therefore it exhibits a out-of-band quasiseparable structure.  Moreover, after the zeroing the second row of $R$ also reveals an out-of-band quasiseparable representation. 
  At the  end of the first step the  matrix $R$ has the following profile:
 \vspace{-0.1cm} 
 \[
R=\left[\begin{array}{cccc}
 0& x & 0 & 0 \\ 0& x& x& x \\ & & x& x\\& & & x  \\ \hline  x& x& x& x\\& x \\ & & x\\ & & & x      \end{array} \right] .
 \]
 The quasiseparability rank is maintained  under the zeroing process  and  thus,  if the algorithm concludes successfully at the very end the matrix $R$ looks like 
 \[
R=\left[\begin{array}{cccc}
 0& 0 & 0 & 0 \\ 0& 0& 0& 0\\ 0&0 & 0& 0\\0&0 &0 & 0 \\ \hline  x& x& x& x\\& x& x& x \\ & & x& x \\ & & & x      \end{array} \right] 
 \]
 where the triangular factor in the lower part of $R$ has quasiseparability rank at most $1$.

 The generalization of this approach for quasiseparable matrices  of any order is considered  below in more details. Let $L_0^T\in \mathbb R^{n\times n}$ be an upper quasiseprable matrix of order $\kappa$ with upper generators 
 $\B g_i\in \mathbb R^{\kappa}$, $1\leq i\leq n-1$, $\B h_j\in \mathbb R^{\kappa}$, $2\leq k\leq n$,
$B_\ell\in \mathbb R^{\kappa \times \kappa}$, $2\leq \ell\leq n-1$.  Set $R_0=\left[\begin{array}{c|c}
L_0 & \sqrt{|\alpha|} \end{array} 
\right]^T$ and $\Omega_0$  such that $A=L_0\Omega_0 L_0^T$.  Recall that 
\[
A_\alpha=A+\alpha I_n=R_0 \widehat \Omega_0 R_0^T, \quad \widehat \Omega_0=\diag[\Omega_0, \pm I].
\]
At the first step we consider the matrices 
\[
G_1=\left[\begin{array}{cc}\widehat X_1 B_{\kappa+1}\\\B g_{\kappa+1}^T\end{array}\right] \in \mathbb R^{(\kappa +1)\times \kappa}, 
\quad \widehat X_1=\left[\begin{array}{cc}\B g_1^T B_2 \cdots B_{\kappa}\\\vdots \\\B g_{\kappa}^T \end{array}\right].
\]
Then  we determine a generalized QL factorization of $G_1$ of the form 
\[
H_1 G_1=\left[\begin{array}{cc}\B 0^T\\X_1\end{array}\right], \quad X_1\in \mathbb R^{\kappa \times \kappa}.
\]
The matrix $H_1$ is constructed as product of  plane rotations determined in accordance with the signature of the matrix $\widehat \Omega_0$ in such a way to obtain 
\[
H_1 \widehat \Omega_0(1:\kappa+1, 1:\kappa+1) H_1^T=\Omega_{1/2}
\]
for a suitable  diagonal signature matrix $\Omega_{1/2}$.  The  calculation  of $X_0$, $G_1$, $H_1$ and $X_1$  does not depend on the shift $\alpha$ and, therefore it can be accomplished offline and stored.  It follows that 
\[R_0^T \Omega_0 R_0=R_{1/2}^T \widehat \Omega_{1/2} R_{1/2}
\]
with 
\[
\widehat \Omega_{1/2}=\diag[\Omega_{1/2}, \widehat \Omega_0(\kappa+2\colon n, \kappa+2\colon n)]
\]
and 
\[R_{1/2}=\left[\begin{array}{c|c}
\widehat L_{1/2} & \sqrt{|\alpha|} \end{array} 
\right]^T, \quad 
\]
where $\widehat L_{1/2}^T$ is  generated from $L_0$  by modifying its first $\kappa+1$ rows so that 
\[
\widehat L_{1/2}^T(1:\kappa+1,  \colon)=H_1 R_0(1:\kappa+1,  \colon)
\]
and 
\[
\widehat L_{1/2}^T(1:\kappa+1,  \colon)=\left[\begin{array}{c|c} F_{1/2} & \begin{array}{ccc} \B 0^T \\ \hline \\ [-0.1cm]\begin{array}{cccc} X_1\B h_{\kappa+2}  & X_1 B_{\kappa+2}\B h_{\kappa+3} & \ldots & \end{array}
\end{array} \end{array}
\right].
\]
The matrix $F_{1/2}\in \mathbb R^{(\kappa+1) \times (\kappa +1)}$ incorporates the  buldge. 
Now, according to the given signatures we apply  Givens/hyperbolic plane rotations in the  planes of coordinates $(1, n+1), (2, n+1), \ldots (\kappa+1, n+1)$  in order to  annihilate the first column of $F_{1/2}$.  Notice that 
\[
R_{1/2}([1:\kappa+1,n+1],   \colon)=\left[\begin{array}{c|c} \begin{array}{ccc} \ \ F_{1/2}\\ \\
\sqrt{|\alpha|} & \B 0^T\end{array}& \begin{array}{ccc} \B 0^T \\ \hline \\ [-0.1cm] \left[\begin{array}{cc}X_1\\\B 0^T\end{array}\right] \left[\begin{array}{cccc}
\B h_{\kappa+2}  & B_{\kappa+2}\B h_{\kappa+3} & \ldots & \end{array} \right]
\end{array} \end{array}
\right].
\]
The application of the plane rotations  converts  $R_{1/2}$ into $R_{1}$  with 
\[
R_{1}([1:\kappa+1,n+1],   \colon)=\left[\begin{array}{c|c} \begin{array}{ccc} \begin{array}{cc}\B 0 & F_{1}\end{array}\\ \\
\B s_1^T\end{array}& \begin{array}{ccc} \B 0^T \\ \hline \\ [-0.1cm] \left[\begin{array}{cc}\widehat X_2\\\B z_1^T\end{array}\right] \left[\begin{array}{cccc}
\B h_{\kappa+2}  & B_{\kappa+2}\B h_{\kappa+3} & \ldots & \end{array} \right]
\end{array} \end{array}
\right], 
\]
where $\widehat X_2\in \mathbb R^{\kappa \times \kappa}$ is lower triangular.  Thus, at the end of the first step we find $R_1$ and a  diagonal signature matrix $\widehat \Omega_1$ such that 
\[
R_{1}=\left[\begin{array}{c|c}
\widehat L_{1} & \widehat S_1 \end{array} 
\right]^T, \quad 
\]
with 
\[
\widehat L_{1}^T(1:\kappa+1,  \colon)=\left[\begin{array}{c|c} \begin{array}{c|c} \B 0 & F_1 \end{array} & \begin{array}{ccc} \B 0^T  \\\hline \\ [-0.1cm]\begin{array}{cccc} \widehat X_2\B h_{\kappa+2}  & \widehat X_2 B_{\kappa+2}\B h_{\kappa+3} & \ldots & \end{array}\end{array} \end{array}
\right], 
\]
\[
\widehat S_1^T(1,   \colon)=\left[ \B s_1^T , \B z_1^T\B h_{\kappa+2},  \B z_1^T  B_{\kappa+2}\B h_{\kappa+3},  \ldots, \B z_1^T  B_{\kappa+2} \cdots B_{n-1}\B h_{n}\right]
\]
\[
\widehat S_1^T(\ell,   \colon)= \sqrt{|\alpha|} \B e_\ell^T, \quad 2\leq \ell \leq n, 
\]
and 
\[R_0^T \widehat \Omega_0 R_0=R_{1}^T \widehat \Omega_{1} R_{1}.
\]
At the second step the  reduction procedure is continued by factoring the matrix 
\[
G_2=\left[\begin{array}{cc}\widehat X_2 B_{\kappa+2}\\\B g_{\kappa+2}^T\end{array}\right]
\]
as 
\[
H_2 G_2=\left[\begin{array}{cc}\B 0^T\\X_2\end{array}\right], \quad X_2\in \mathbb R^{\kappa \times \kappa}.
\]
Multiplying $R_1(2:\kappa+2, \colon)$ on the left by $H_2$ has the effect of creating zeros in the second row of 
$R_{1}$ by moving the bulge one position  down  along the main  diagonal.  
Notice that 
\[
\widehat L_{3/2}^T(1:\kappa+2,  2\colon n)=\left[\begin{array}{c|c} F_{3/2} & \begin{array}{cccc} \B 0^T \\ \B 0^T \\\hline \\ [-0.1cm]\begin{array}{cccc} X_2\B h_{\kappa+3}  & X_2 B_{\kappa +3}\B h_{\kappa+4} & \ldots & \end{array}
\end{array} \end{array}
\right]
\]
with $F_{3/2}\in \mathbb R^{(\kappa +2)\times (\kappa+1)}$. 
The first column $F_{3/2}$  can then  be zeroed by  transformations in the planes of coordinates $(1, n+2), (2, n+2), \ldots (\kappa+2, n+2)$. As the process goes on, the matrix $F_{\ell-1/2}$ continues to grow one row at a time
in the  first $\kappa +1$  steps after which it begins to move down the main diagonal  with constant size $(2\kappa+1)\times (\kappa+1)$.
Intuitively speaking,  at the regime the generic step of the zeroing process looks as follows:
\begin{enumerate}
\item\label{r1} Compute the generalized QL decomposition of $G_\ell=\left[\begin{array}{cc}\widehat X_\ell B_{\kappa+\ell}\\\B g_{\kappa+\ell}^T\end{array}\right]$ 
as  
$H_\ell G_\ell=\left[\begin{array}{cc}\B 0^T\\X_\ell\end{array}\right]$, with $ X_\ell\in \mathbb R^{\kappa \times \kappa}$, and  $\widehat H_\ell \widehat \Omega_{\ell-1} \widehat H_\ell^T =\widehat \Omega_{\ell-1/2}$, $\widehat H_\ell=I_{\ell-1}\oplus H_\ell\oplus I_{2n-\ell-\kappa}$; 
\item\label{r2}  Compute $R_{\ell-1/2}=\widehat H_\ell R_{\ell-1}$; 
\item\label{r3} For $k=\ell-\kappa\colon \ell$  perform 
\[[R_{\ell-1/2}, \widehat \Omega_{\ell-1/2}]=ZEROING(R_{\ell-1/2}, \widehat \Omega_{\ell-1/2}, k, \ell).
 \]
\item\label{r4}  For $k=\ell+1:\ell+\kappa$ determine 
 \[[R_{\ell-1/2}, \widehat \Omega_{\ell-1/2}]=ZEROING(R_{\ell-1/2}, \widehat \Omega_{\ell-1/2}, k, \ell)
 \]
 by calculating $\widehat X_{\ell+1}$.
 \item\label{r5} Set $\widehat \Omega_{\ell}=\widehat \Omega_{\ell-1/2}$ and $R_{\ell}=R_{\ell-1/2}$.
\end{enumerate}

The zeroing scheme performed at steps \ref{r3} and \ref{r4}  acts on the submatrix 
$R_{\ell-1/2}([ \ell-\kappa:\ell+\kappa,n+\ell],   \colon)$ given by 
\[
\left[\begin{array}{c|c} \begin{array}{ccc} \ \ F_{\ell-1/2}\\ \\
\sqrt{|\alpha|} & \B 0^T\end{array}& \begin{array}{ccc} 0_{\kappa+1, n-\kappa-\ell+1} \\ \hline \\ [-0.1cm] \left[\begin{array}{cc}X_\ell\\\B 0^T\end{array}\right] \left[\begin{array}{cccc}
\B h_{\kappa+\ell+1}  & B_{\kappa+\ell+1}\B h_{\kappa+\ell+2} & \ldots & \end{array} \right]
\end{array} \end{array}
\right]
\]
where $X_\ell\in \mathbb R^{\kappa\times \kappa}$  is lower triangular and, moreover, $F_{\ell-1/2}\in \mathbb R^{(2\kappa+1)\times (\kappa+1)}$  is  a trapezoidal matrix such that 
\[
F_{\ell-1/2}=\left[\begin{array}{ccc} *\\ \vdots & \ddots \\ \raisebox{-3pt}{*} & \raisebox{0pt}{\ldots} & \raisebox{-3pt}{*} \\   \vdots & \vdots & \vdots \end{array}\right]. 
\]
The application of the plane rotations at steps \ref{r3} and \ref{r4}  modifies this submatrix  so that 
\[
R_{\ell}([\ell-\kappa:\ell+\kappa,n+\ell],   \colon)=\left[\begin{array}{c|c} \begin{array}{ccc} \begin{array}{cc}\B 0 & F_{\ell}\end{array}\\ \\
\B s_\ell^T\end{array}& \begin{array}{ccc} 0_{\kappa+1, n-\kappa-\ell+1}\\ \hline \\ [-0.1cm] \left[\begin{array}{cc}\widehat X_{\ell+1}\\\B z_\ell^T\end{array}\right] \left[\begin{array}{cccc}
\B h_{\kappa+\ell+1}  & B_{\kappa+\ell+1}\B h_{\kappa+\ell+2} & \ldots & \end{array} \right]
\end{array} \end{array}
\right], 
\]
where $\widehat X_{\ell+1}$ is lower triangular. 
In this way  after $n-1$ steps at the very end  of the zeroing scheme  we find a matrix $R_{n-1}$ and a signature matrix $\widehat \Omega_{n-1}$ such that 
\[
A_\alpha=A+\alpha I_n=R_0^T \widehat \Omega_0 R_0=R_{n-1}^T \widehat \Omega_{n-1} R_{n-1}
\]
where 
\[
R_{n-1}=\left[\begin{array}{c|c}
0 & \widehat S_{n-1} \end{array} 
\right]^T, \quad 
\]
and $\widehat S_{n-1}=L_\alpha $ is upper triangular with a suitable  rank structure. Specifically,  we find that 
\begin{theorem}\label{main}
Whenever the zeroing scheme \ref{r1}--\ref{r5} goes to completion,  it computes 
vectors $\B s_i$, $1\leq i\leq n$, $\B z_i$,  $1\leq i\leq n-\kappa-1$, and a signature matrix $\Omega_\alpha$ such that $A_\alpha=L_\alpha \Omega_\alpha L_\alpha^T$, $L_\alpha^T$ is upper triangular and its entries are given by 
\[
(L_\alpha^T)_{i,j}=\left\{\begin{array}{lll}
0 \ {\rm if} \ i>j; \\
(\B s_i)_{j-i+1} \ {\rm if} \ i\leq j\leq i+\kappa; \\
\B z_i^T  B_{\kappa +i +1}\cdots B_{j-1} \B h_j \ {\rm if}  \ i+\kappa<j\leq n.
\end{array}\right.
\]
\end{theorem}
Therefore, $L_\alpha^T$  can be represented by  a band profile  of bandwidth $\kappa$ and  a quasiseparable
structure of order $\kappa$ above the band. Such matrices are referred as out-of-band quasiseparable matrices in \cite{EIDELMAN2008266}.
The computational  cost  of the  resulting updating procedure  can be estimated as $O(n (\kappa^2 + \theta(\kappa) +\vartheta(\kappa))$, where $\theta(k)$  and $\vartheta(\kappa)$ denote  the cost of performing the matrix multiplication $\widehat X_\ell B_{\kappa +\ell}$ and  of factoring  the matrix 
$G_\ell=\left[\begin{array}{cc}\widehat X_\ell B_{\kappa+\ell}\\\B g_{\kappa+\ell}^T\end{array}\right]$
in its QL decomposition, respectively. Therefore, the worst-case complexity is $O(n\kappa^3)$. 
However, recall  that $\widehat X_\ell$ is lower triangular of size $\kappa$. Therefore, if $B_\ell=I_\kappa$, $2\leq \ell\leq n-1$  we find that the overall cost is $O(n\kappa^2)$ flops which is optimal both respect to  the size $n$ and  the quasiseparability order $\kappa$.  
Finally,it could be shown that $S_{n-1}$ is  indeed quasiseparable of order $\kappa$ and its minimal representation can be generated from the out-of-band structure in some post-processing stage. 

\section{Numerical Results}\label{four}
Sequences of shifted  rank-structured linear systems are encountered in the solution of time evolution problems  for linear operator  equipped with nonlocal boundary conditions.   Here we consider  two different model problems. 

The first  example is Fredholm-type time-dependent linear integro-differential equation  given by:
 \begin{eqnarray}\label{prob1}
&&\frac{\partial u(x,t)}{\partial t} =\alpha u(x,t) + \int_0^1 B(x, s)u(s,t)\,\mathrm{d}s  =\mathcal L u(x,t) \quad 0\leq x,t \leq 1,   
\end{eqnarray}
equipped with the integral boundary condition
\begin{eqnarray}\label{prob2}
  &&\displaystyle \int_0^1 u(x,t)\,\mathrm{d}t = g(x).
\end{eqnarray}
Here 
\[
\mathcal L w(x)=\alpha w(x) + \int_0^1 B(x, s)w(s)\,\mathrm{d}s,  
\]
is a linear operator acting in some Banach space. 
Assume that the kernel  $B(x, s)$  is semiseparable, 
that is, 
\[
B(x, s)=\left\{ \begin{array}{ll}\sum_{i=1}^\kappa \phi_i(x)\psi_i(s)  \ {\rm if}\ x\leq s; \\ \\
\sum_{i=1}^\kappa \phi_i(s)\psi_i(x)  \ {\rm if}\ x\geq  s.
\end{array}\right.
\]
The problem \eqref{prob1}-\eqref{prob2} can be  partially discretized by  the method of lines using a semidiscretization in space   with equispaced points $x_i=(i-1)h$, $1\leq i\leq n$, $h=1/(n-1)$. Let us denote $u_i(t)=u(x_i, t)$ , $\B u(t)=\left[ u_1(t), \ldots, u_n(t)\right]^T$ and consider the approximation of the integral in \eqref{prob1} by means of the trapezoidal rule.   In this way we obtain the semidiscretized formulation 
\begin{eqnarray}\label{prob3}
&&\frac{d \B u(t)}{d t} =L \B u(t) \quad 0\leq t \leq 1,   \quad L=(l_{i,j})\in \mathbb R^{n\times n}
\end{eqnarray}
equipped with the integral boundary condition
\begin{eqnarray}\label{prob4}
  &&\displaystyle \int_0^1 \B u(t)\,\mathrm{d}t = \left[g(x_1), \ldots,  g(x_n)\right]^T.
\end{eqnarray}
The  matrix $L$ satisfies $L=A D$, where $D=\diag[1/2, 1, \ldots, 1, 1/2]$ and $A=(a_{i,j})$ with 
$a_{i,j}= h B(x_i, x_j)$ for $i\neq j$.  Therefore $A$ is symmetric  and  moreover \eqref{prob3}-\eqref{prob4} can be rewritten as 
\begin{eqnarray}\label{prob5}
&&\frac{d \B {\widehat u}(t)}{d t} =\widehat L  \B {\widehat u}(t) \quad 0\leq t \leq 1, 
\end{eqnarray}
equipped with the integral boundary condition
\begin{eqnarray}\label{prob6}
  &&\displaystyle \int_0^1 \B {\widehat u}(t)\,\mathrm{d}t = \left[\hat g(x_1), \ldots,  \hat g(x_n)\right]^T =\B {\widehat g}, 
\end{eqnarray}
where $\B {\widehat u}(t)=D^{1/2} \B u(t)$, $\widehat L=D^{1/2} A D^{1/2}$ and $\hat g(x_i)= \sqrt{D_{i,i}} g(x_i)$, $1\leq i\leq n$.

In \cite{Boito1} it is proved that  the nonlocal boundary value problem \eqref{prob5}-\eqref{prob6}
admits as unique solution
\begin{equation}\label{l3}
 \mathbf{ \widehat u}(t)= q(t,\widehat L) \B {\widehat g},   \quad q(t,w)= \frac{w e^{w t}}{e^w -1}, \qquad
 w \in \mathbb{C} \setminus \{0, \pm 2\pi \mc, \pm 4 \pi \mc, \ldots\}
 \end{equation}
A family of rational approximations of the solution  $\B {\widehat u}(t)$  of \eqref{prob5}-\eqref{prob6} has been devised and studied in the papers \cite{Aceto2025ER,Boito,Boito1,Boito2025}.  Any approximant 
 $\mathbf u_{\ell, N}(t)$ can be described as follows: 
\[
    \mathbf {\widehat u}_{\ell, N}(t)=\mathbf p_\ell(t)+\mathbf s_{\ell, N}(t)
\]
with 
\[
\mathbf p_\ell(t)=\sum_{k=0}^{2\ell+1} \frac{B_k(t)}{k!}  
\B {\widehat g_k}, \quad \B {\widehat g_k}=\widehat L^k \B {\widehat g}
\]
and
\[
\mathbf s_{\ell, N}(t)=(-1)^{\ell}2\sum_{k=1}^{N}\left(\frac{1}{2 \pi k}\right)^{2\ell}
\left (\Sigma_k(\mathbf {\widehat g_{2\ell +1}})\cos(2 \pi k t)+ \Upsilon_k(\mathbf {\widehat g_{2\ell +1}})\sin(2 \pi k t)\right )
\]
where $B_m(t)$ are the well-known Bernoulli polynomials:
$$
B_m(t)=\sum_{j=0}^m\frac1{j+1}\sum_{k=0}^j(-1)^k
\left(\begin{array}{c}j\\k\end{array}\right)(t+k)^m.
$$
and 
\begin{equation}\label{shiftsys}
\left\{ \begin{array}{ll}
\Sigma_k(\mathbf {\widehat g_{2\ell +1}})=\widehat L (\widehat  L^2+ (2 \pi k)^2  I )^{-1}\mathbf {\widehat g_{2\ell +1}},\\
\Upsilon_k(\mathbf {\widehat g_{2\ell +1}})=(2 \pi k)^{-1} \widehat  L \Sigma_k(\mathbf {\widehat g_{2\ell +1}}), 
\end{array}\right. \quad k=1,2,\ldots, N.
\end{equation}
The computational bulk is given by the solution of the shifted linear systems \eqref{shiftsys} with coefficient matrix $\widehat L^2$.  An initial Cholesky-like decomposition  of $\widehat L^2$ can be derived from the  QR factorization of $\widehat L$.  A structured QR  algorithm working on generators is presented in \cite{EIDELMAN2002419}. We have implemented this algorithm together with our updating procedure in \texttt{MATLAB} and then we have performed  several numerical experiments to demonstrate  the efficiency and the accuracy of the proposed method for solving \eqref{prob5}-\eqref{prob6}.
In our numerical experiments we consider as reference  solution of \eqref{prob5}-\eqref{prob6}
\[
\B {\widehat u}(t)=({\tt expm}(\widehat L) -{\tt eye}(s))\textbackslash ({\tt expm}(t \widehat L) \widehat L\B {\widehat g})
\]
computed  from \eqref{l3} using the \texttt{MATLAB} operator "backslash" and the built-in function ${\tt expm}$.  Problem data are specified  as follows: 
\[
\phi_k(x)=\sin(2\pi k x), \ \psi_k(x)=\cos((7-2k)\pi x), \ k=1, 2, 3; \quad g(x)=(x-1)/x(+1).
\]
Therefore, $\widehat L$ is quasiseparable of  quasiseparability order 3 represented in semiseparable form. Our implementation of the structured QR algorithm  devised in \cite{EIDELMAN2002419} returns an initial Cholesky  factorization of the matrix $\widehat L^2=L_0 L_0^T$ with $L_0$ lower quasiseparable of quasiseparability order 6.

\begin{figure}
\centering
 \includegraphics[width=0.5\textwidth]{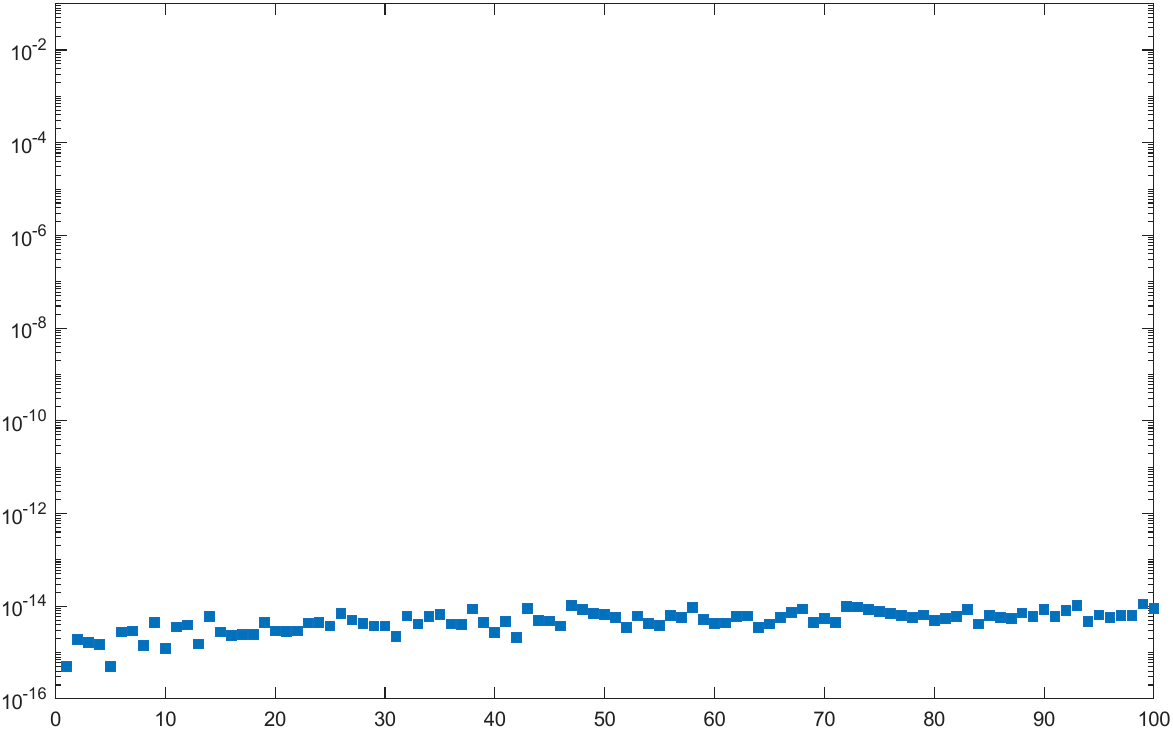}
 \caption{Plot of the backward error in the solution of the shifted linear systems \eqref{shiftsys}  for $n=1024$ and $1\leq k\leq 100$.}
\label{f1}
\end{figure}

\begin{figure}
 \subfloat[]{\includegraphics[width=0.5\textwidth]{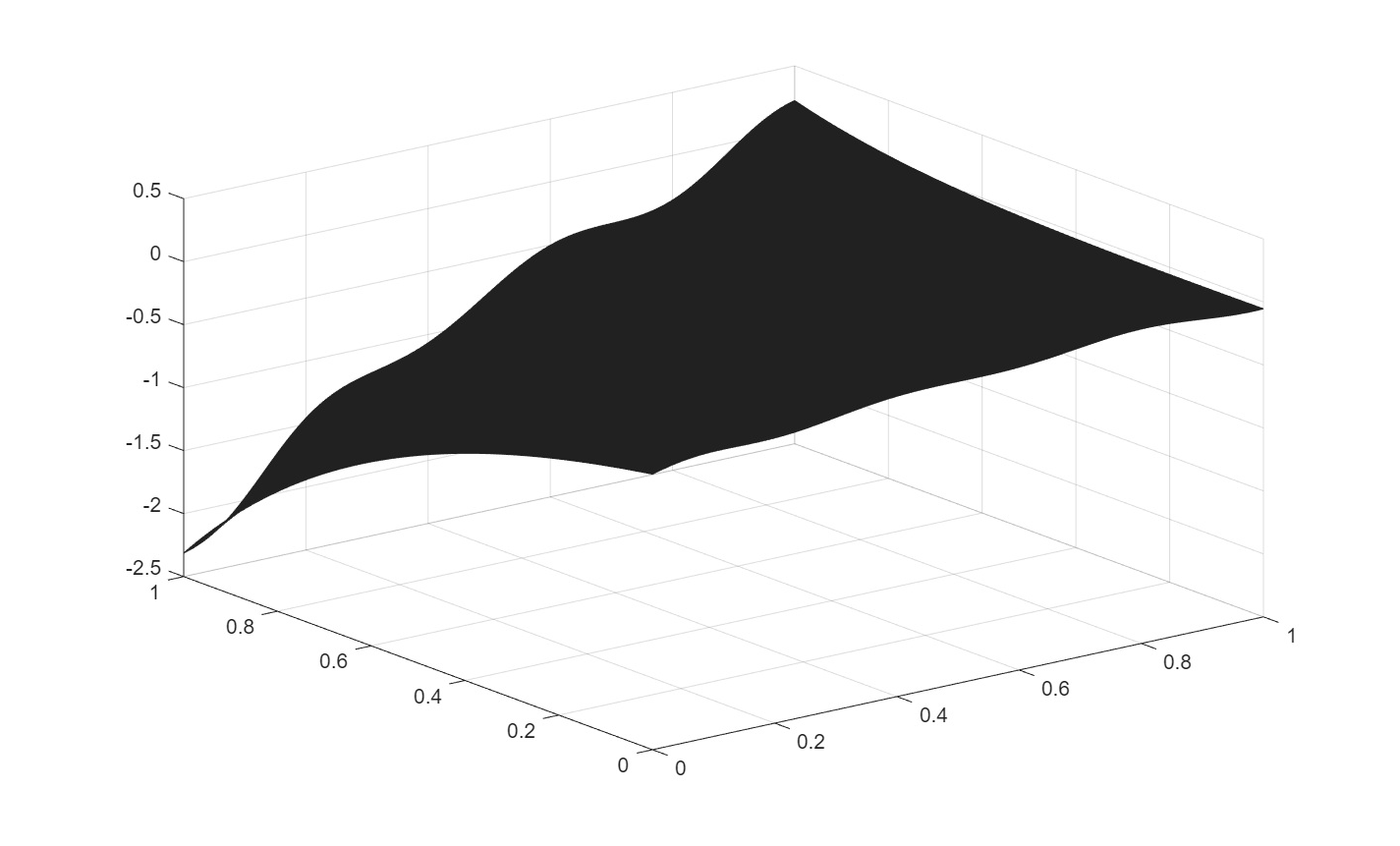}}
 \hfill
 \subfloat[]{\includegraphics[width=0.5\textwidth]{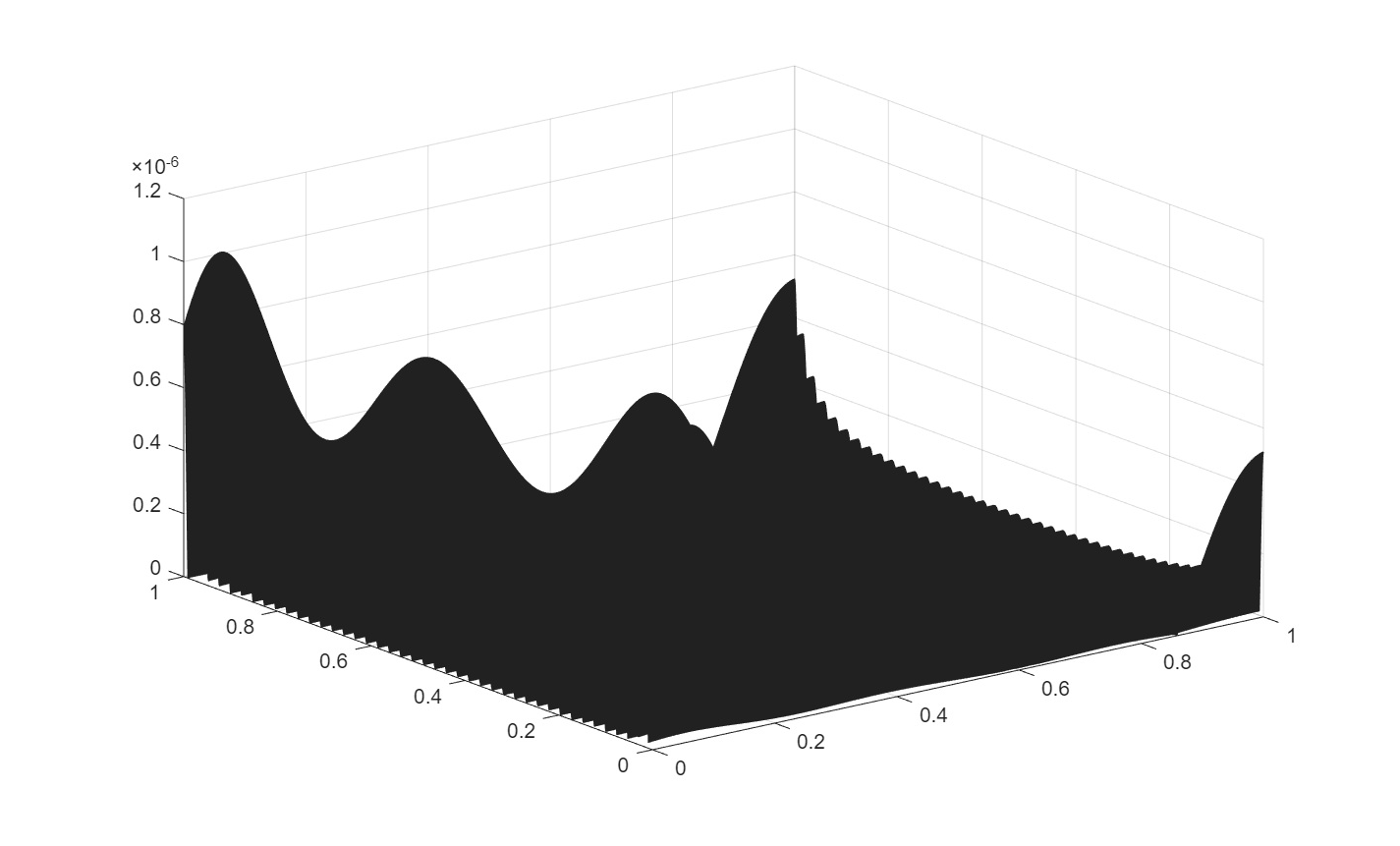}}
 \caption{Plots of the  approximated solution and absolute error for  $n=1024$ with $\ell=1$ and 
 $N=20$.}
\label{f2}
\end{figure}

\begin{figure}
 \subfloat[]{\includegraphics[width=0.5\textwidth]{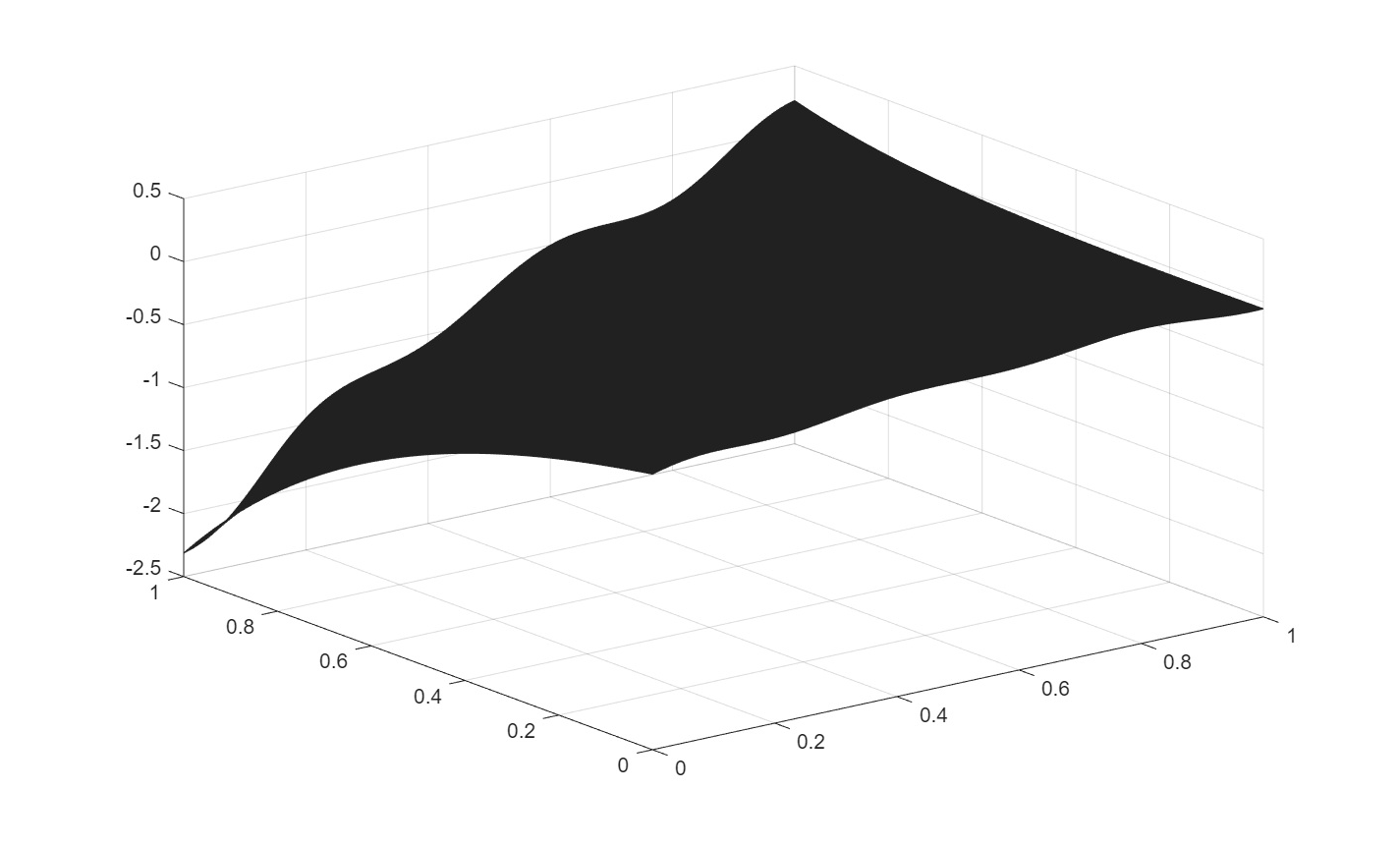}}
 \hfill
 \subfloat[]{\includegraphics[width=0.5\textwidth]{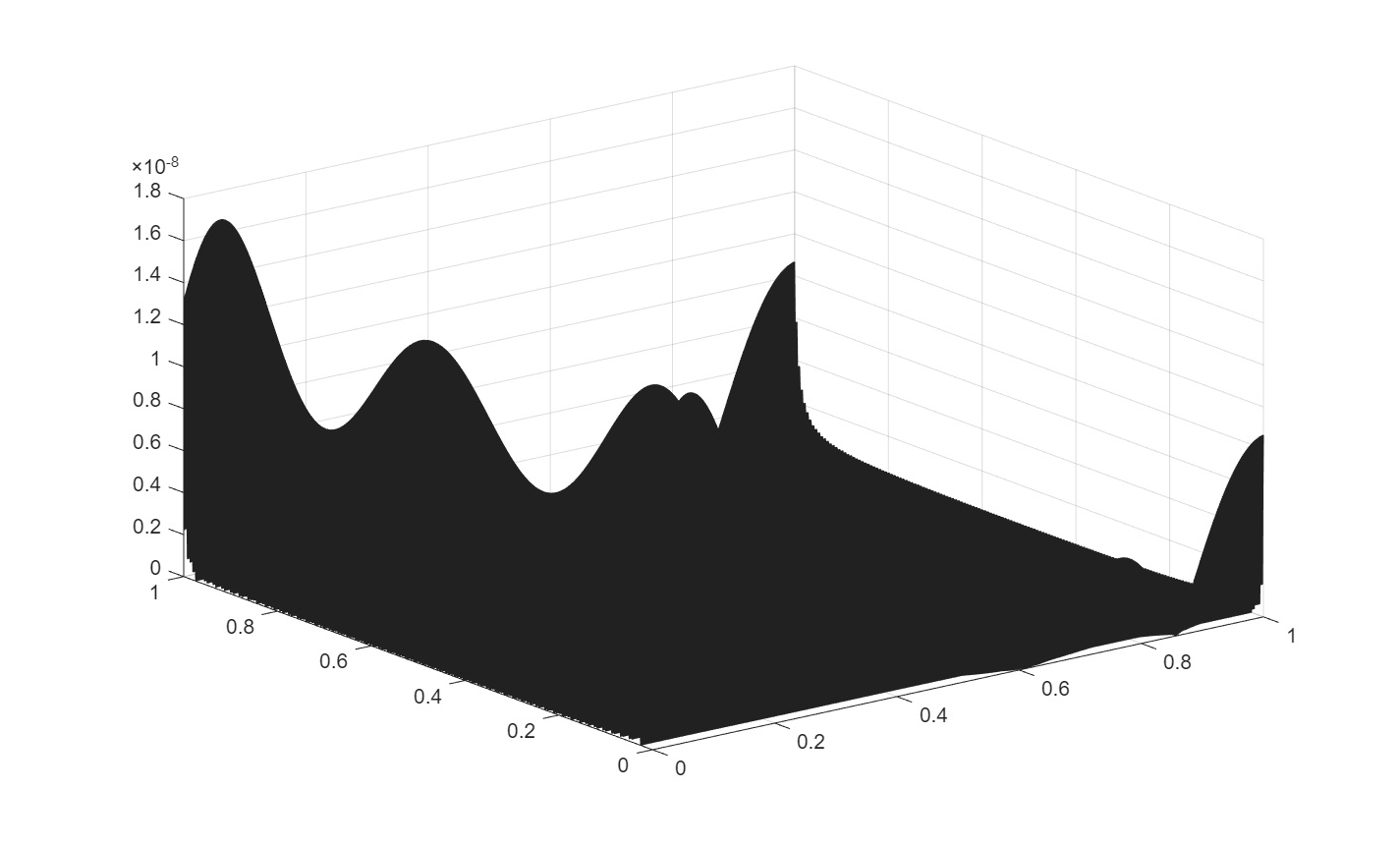}}
 \caption{Plots of the  approximated solution and absolute error for  $n=1024$ with $\ell=1$ and 
 $N=80$.}
\label{f3}
\end{figure}

First, we have checked  the accuracy of our proposed numerical method for solving  the shifted linear systems \eqref{shiftsys}.  In Figure \ref{f1} we show  the backward error 
\[
err=\frac{\parallel \widehat  L^2+ (2 \pi k)^2  I  -S_{n-1}\widehat \Omega(n+1\colon 2n, n+1\colon 2n) S_{n-1}^T\parallel_2}{\parallel \widehat  L^2+ (2 \pi k)^2  I \parallel_2 }
\]
for $n=1024$ and $k=1, \ldots, 100$.  Numerical results clearly indicate the backward stability of the resulting algorithm. Notice that in the considered  situation since $\widehat L^2$ is positive definite  and the shifts are positive only Givens plane rotations have been applied. 
Then, we have tested the accuracy of the approximation  $\mathbf {\widehat u}_{\ell, N}(t)$ of the reference solution $\B {\widehat u}(t)$  of \eqref{prob5}-\eqref{prob6}. In order to illustrate the approximation error, we have plotted   the function $(\mathbf {\widehat u}_{\ell, N}(t_j))_i$ and the absolute error $|(\mathbf {\widehat u}_{\ell, N}(t_j))_i-(\B {\widehat u}(t_j))_i|$ over the grid $X\times T$ where  $X=\{x_1, \ldots, x_n\}$ and 
$T=\{t_1, \ldots, t_n\}$ with $t_i=(i-1)h$. In Figure \ref{f2} and \ref{f3} we show the results for $n=1024$, $\ell=1$, and $N=20$ and $N=80$, respectively. 
We see that the approximation is already effective for $N=20$.  The error  typically 
increases at the endpoints of the time interval due to the  possible occurrence of points of discontinuity since the solution of the continuous problem  \eqref{prob1}-\eqref{prob2} is not generally  periodic.  The reduction in the approximation error for $N=80$ is in accordance with the convergence  estimates provided in  \cite{Aceto2025ER}. 
The approximation in a   region far-away  from the endpoints of the time interval might be further improved by means of the acceleration technique presented in \cite{Aceto2025ER}.

The second  example is a classical  diffusion problem of heat in a square region $\Omega=[0,1]\times [0,1]$ with an integral condition.  The problem looks like: 
\begin{eqnarray*}\label{heatproblem}
    &&\frac{\partial u(x,y,t)}{\partial t} = \alpha^2 \left(\frac{\partial^2 u(x,y,t)}{\partial x^2} + \frac{\partial^2 u(x,y,t)}{\partial y^2}\right) + f(x, y) u(x,y, t), \ (x,y)\in \Omega  \\
  &&\displaystyle \int_0^1 u(x, y, t)\,\mathrm{d}t = g(x, y),  \label{hl2}
\end{eqnarray*}
with homogeneous  Dirichlet boundary conditions
\begin{equation}\label{heatinitialdata}
u(0, y, t)=u(1, y,t)=u(x,0,t)=u(x,1,t)=0, \quad 0\leq x,y\leq 1, t\geq 0.
\end{equation}
Here 
\[
\mathcal L w(x,y)=\alpha^2 \left(\displaystyle\frac{\partial^2 w(x,y)}{\partial x^2} + \frac{\partial^2 w(x,y)}{\partial y^2}\right) + f(x,y) w(x,y)
\]
where $w(x,y)$ satisfies  the 
boundary conditions $w(0,y)=w(1,y)=w(x,0)=w(x,1)=0$.

\begin{figure}
\begin{minipage}[t]{.45\linewidth}
\vspace{0pt}
\centering
\includegraphics[width=1.5in]{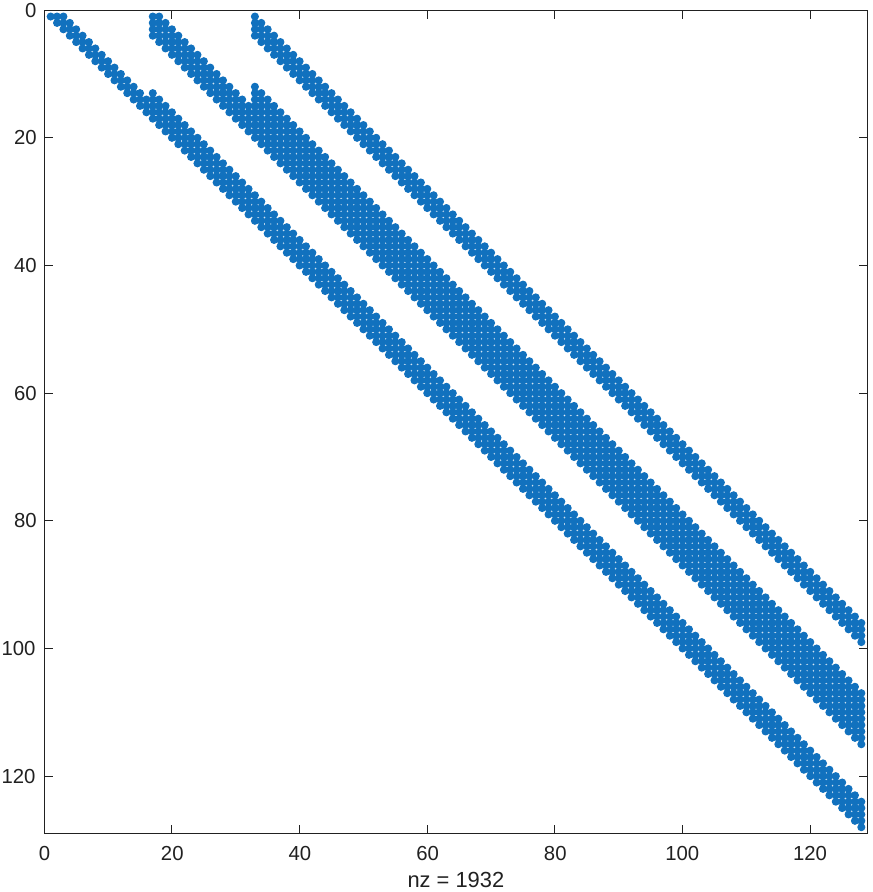}
\end{minipage}%
\begin{minipage}[t]{.25\linewidth}
\vspace{0pt}
\centering
\begin{tabular}{|c|c|c|c|c|c|}
\hhline{-|-----}
\multicolumn{1}{|c|}{\backslashbox{$\beta$}{$n$}}
&8 & 16 & 32& 64 & 128\\
\hhline{-|-----}
$(2\pi)^2$ & 1 & 1 & 1& 1& 1\\
\hhline{-|-----}
$(4\pi)^2$ & 1 & 1 & 1& 1& 1\\
\hhline{-|-----}
$(6\pi)^2$ & 1 & 1 & 1& 1& 1\\
\hhline{-|-----}
$(8\pi)^2$ & 2& 1 & 1& 1& 1\\
\hhline{-|-----}
$(10\pi)^2$ & 2 & 2 & 1& 1& 1\\
\hhline{-|-----}
\end{tabular}
\end{minipage}
\caption{ The spy plot on the left shows the sparsity pattern of the matrix $R_0$.  The table on the right gives the number of iterations required by {\tt pcg} of MATLAB  with prescribed tolerance  \num{1.0e-8}}\label{ffp}
\end{figure}

 By using a finite difference discretization over the uniform grid  $\{(x_i, y_j)\in \mathbb R^2\colon x_i=y_i=ih, 0\leq i\leq n+1\}$, $h=1/(n+1)$,  of  the unit square we may write the discretized version of the problem in the form \eqref{prob5}  with 
 \[
 \widehat L=  \displaystyle\left(\frac{\alpha}{h}\right)^{2} T +\diag\left[f(x_1, y_1), \ldots, f(x_1, y_n), f(x_2, y_1), \ldots,  f(x_n, y_n)\right]
 \]
 and  $T\in \mathbb R^{N\times N}$, $N=n^2$,  is the classical approximation of the $2D-$Laplacian, i.e., 
\[
T=\left[\begin{array}{cccc}
T_1& I_N\\ I_N & T_1& \ddots \\ & \ddots & \ddots & I_N\\ & & I_N & T_1
\end{array}\right],  \quad  T_1=\left[\begin{array}{cccc}
-4& 1\\ 1 & -4& \ddots \\ & \ddots & \ddots & 1\\ & & 1 & -4
\end{array}\right] \in \mathbb R^{n\times n}.
\]
The efficient solution of the shifted linear systems \eqref{shiftsys} is performed as follows.  We first compute  the Cholesky factorization of the matrix $A=\widehat L^2$, say $A=\widehat R_0 \widehat R_0^T$.  Then we determine a sparse approximation $R_0$ of $\widehat R_0$ and  for any shift $\beta>0$ we calculate the incomplete  Cholesky-factorization of $R_0  R_0^T  + \beta I$   according to the sparsity pattern of $R_0$ by using Procedure {\textsc{INUPDATE}} at a linear cost.   Finally, the  incomplete triangular factor  is used as preconditioner in the  MATLAB function {\tt pcg} which implements the preconditioned conjugate gradient algorithm for solving the linear system with coefficient matrix $A+\beta I$.  By neglecting the amount of work needed to find the initial factorization,  we can estimate the overall computational cost as $O(N \cdot it)$,   where $it$ is the number of iterations of {\tt pcg} required to achieve the prescribed tolerance. In Figure \ref{ffp} we  illustrate the sparsity pattern of the matrix $R_0$ and the number of iterations  carried out by the preconditioned conjugate gradient method applied to  the shifted systems generated with $f(x,y)=\sin(x)\cos(y)$ and $\alpha^2=\num{1.18e-6}$ --the thermal diffusivity of sand--  for $n=2^j$, $3\leq j\leq 7$. 

It is clearly seen from Figure \ref{ffp} that our proposed method solves   each shifted linear system  with linear  complexity w.r.t. the size of the coefficient matrix.

\section{Conclusions and Future Work}\label{five}
A novel method for efficiently updating the Cholesky-like factorization of a scalar quasiseparable matrix has been introduced . The method can be applied for the solution of sequences of shifted linear systems occurring in  some numerical schemes for boundary value problems  equipped with nonlocal boundary conditions.  The results of numerical experiments are reported to confirm the robustness and the effectiveness of the proposed method. Future work is concerned with the extension of the updating procedure to the block case and the use of the Cholesky-like factorization approach as the workhorse for the design of a fast and reliable bisection algorithm for eigenvalue computation of quasiseparable matrices.

\begin{acknowledgements}
The author is member of INDAM/GNCS.
\end{acknowledgements}

\bibliographystyle{plain} 
\bibliography{bibliography}
\end{document}